\date{}
\documentclass [11pt]{article}
\usepackage{amsfonts,amsmath,graphicx}

\title{\vspace{-0.3cm}Triangle packings and $1$-factors in oriented graphs}
\author{
Peter Keevash \thanks{School of Mathematical Sciences,
Queen Mary, University of London, Mile End Road, London E1 4NS, UK.
Email: p.keevash@qmul.ac.uk.
Research supported in part by NSF grant DMS-0555755.}
\and
Benny Sudakov \thanks{Department of Mathematics, UCLA, Los Angeles, 90095. Email:
bsudakov@math.ucla.edu.
Research supported in part by NSF CAREER award DMS-0546523,
and a USA-Israeli BSF grant.}
}

\newtheorem{theo}{Theorem}[section] 

\newtheorem{lemma}[theo]{Lemma}
\newtheorem{coro}[theo]{Corollary}

\newcommand{\mb}[1]{\mathbb{#1}}
\newcommand{\nib}[1]{\noindent {\bf #1}}

\newcommand{\sm}{\setminus}
\newcommand{\ov}{\overline}
\newcommand{\eps}{\varepsilon}

\newcommand{\sub}{\subseteq}

\def\qed{\hfill $\Box$}

\begin{document}
\maketitle

\begin{abstract}
An oriented graph is a directed graph which can be obtained from a simple
undirected graph by orienting its edges. In this paper we show that
any oriented graph $G$ on $n$ vertices
with minimum indegree and outdegree at least $(1/2-o(1))n$ contains
a packing of cyclic triangles covering all but at most $3$ vertices.
This almost answers a question of Cuckler and Yuster and is best possible,
since for $n \equiv 3$ mod
$18$ there is a tournament with no perfect triangle packing
and with all indegrees and outdegrees
$(n-1)/2$ or $(n-1)/2 \pm 1$. Under the same hypotheses, we also
show that one can embed any prescribed almost $1$-factor, i.e.
for any sequence $n_1,\cdots,n_t$ with $\sum_{i=1}^t n_i \le n-O(1)$
we can find a vertex-disjoint collection of directed cycles with
lengths $n_1,\cdots,n_t$. In addition, under quite general conditions on the $n_i$ we can
remove the $O(1)$ additive error and find a prescribed $1$-factor.
\end{abstract}

\section{Introduction}

A classical result of Extremal Combinatorics, Dirac's Theorem \cite{D}, states that
a graph $G$ on $n \ge 3$ vertices with minimum degree at least $n/2$ contains a
Hamiltonian cycle, i.e. a cycle that passes through every vertex of $G$.
This motivates the general question of determining what minimum degree condition
one needs to find a certain structure in a graph. An important result
of this type is the Hajnal-Szemer\'edi Theorem \cite{HS}, which states that
if $G$ is a graph on $n$ vertices with minimum degree at least $(1-1/r)n$
and $r$ divides $n$ then $G$ has a perfect $K_r$-packing, i.e., a collection
of vertex-disjoint copies of the complete graph $K_r$ on $r$ vertices which covers all
the vertices of $G$. (The case $r=3$ was
obtained earlier by Corr\'adi and Hajnal \cite{CoH}.)
This was generalised to packings of arbitrary graphs
by Koml\'os, Sark\"ozy and Szemer\'edi \cite{KSS2}. Confirming
a conjecture of Alon and Yuster \cite{AY}, they
proved that for every graph $H$ there is a constant $C$ such
that if $G$ is a graph on $n$ vertices with minimum degree at least
$(1-1/\chi(H))n+C$ and $|V(H)|$ divides $n$ then $G$ has a perfect $H$-packing.
(Here $\chi(H)$ denotes the chromatic number of $H$.)
Finally, K\"uhn and Osthus \cite{KO1} determined the minimum degree needed
to find an $H$-packing up to an additive constant: it is $(1-1/\chi^*(H))n+O(1)$,
where $\chi^*(H)$ is a rational number in the range $(\chi(H)-1,\chi(H)]$
that can be calculated when $H$ is given.
Another packing result that is closely related to the topic of this
paper was obtained by Aigner and Brandt \cite{AB} and in a slightly weaker form by Alon and Fischer \cite{AF}.
Verifying a conjecture of Sauer and Spencer \cite{SS},
they proved that a graph $G$ on $n$ vertices with minimum
degree at least $(2n-1)/3$ contains any graph $H$ on $n$ vertices with maximum
degree at most $2$.

It is very natural to ask whether these results have analogues for directed graphs.
Here instead of degree one may consider the {\em minimum semi-degree}
$\delta^0(G)= \min \big(\delta^+(G),\delta^-(G)\big)$,
where $\delta^+(G)$ is the minimum outdegree and $\delta^-(G)$ is the minimum
indegree of a digraph $G$. A directed version of Dirac's theorem was obtained by
Ghouila-Houri \cite{G}, who showed that any digraph $G$ on $n$ vertices with minimum semi-degree
at least $n/2$ contains a Hamilton cycle. (When referring to paths and cycles in directed
graphs we always mean that these are directed, without mentioning this explicitly.)
This result is very closely related to Dirac's theorem, as (some)
extremal digraphs can be obtained from extremal graphs for Dirac's theorem by
replacing each edge by a pair of arcs, one in each direction; the proof of the upper bound
is more complicated, but not unduly so. However, the situation becomes more complicated if one considers
\emph{oriented graphs}. An oriented graph is a directed graph that can be obtained from a (simple)
undirected graph by orienting its edges. The question concerning the analogue of
Dirac's theorem in this case was raised by Thomassen \cite{T1}, who asked what minimum
semi-degree forces a Hamilton cycle in an oriented graph.
Over the years since the question was posed, a series of improving bounds were
obtained in \cite{T2, T3, H, HT}, until an asymptotic solution of $(3/8+o(1))n$ was given by
Kelly, K\"uhn and Osthus \cite{KelKO}. Soon after that an exact answer $\left\lceil \frac{3n-4}{8} \right\rceil$
for $n$ sufficiently large was proved by Keevash, K\"uhn and Osthus \cite{KeeKO}.
This result was recently extended by Kelly, K\"uhn and Osthus \cite{KelKO2}, who showed
that the same minimum semi-degree condition guarantees that $G$ is pancyclic (contains
directed cycles with all lengths $\ell$, $3 \le \ell \le n$).

In this paper we mainly study cycle packings in oriented graphs, although
in the concluding remarks we discuss directed graphs as well.
Our starting point is the following question posed independently by
Cuckler \cite{C} and Yuster \cite{Y}. A tournament is an
orientation of a complete graph. It is regular
if every vertex has equal indegree and outdegree.

\nib{Question.} Does a regular tournament
on $n$ vertices with $n \equiv 3$ mod $6$ have a perfect packing
of cyclic triangles?

We obtain the following general result, which in the case of tournaments
`almost' answers this question.

\begin{theo} \label{triangles}

(i) There is some real $c>0$ so that for sufficiently large $n$,
any oriented graph $G$ on $n$ vertices
with minimum indegree and outdegree at least $(1/2-c)n$ contains
a packing of cyclic triangles covering all but at most $3$ vertices.

(ii) If $n \equiv 3$ mod $18$ then there is a tournament $T$
which does not have a perfect packing of cyclic triangles,
in which every vertex has indegree and outdegree
$(n-1)/2$ or $(n-1)/2 \pm 1$.

\end{theo}

Our second result is an attempt to prove a directed analogue of the results of
Aigner-Brandt and Alon-Fischer. It shows that the same
semidegree condition as above allows one to
cover all but a constant number of vertices by cycles of prescribed lengths.

\begin{theo} \label{almost-factor}
There exist constants $c, C>0$ such that for $n$ sufficiently large,
if $G$ is an oriented graph on $n$ vertices
with minimum indegree and outdegree at least $(1/2-c)n$
and $n_1,\cdots,n_t$ are numbers with $\sum_{i=1}^t n_i \le n-C$
then $G$ contains vertex-disjoint cycles of length $n_1,\cdots,n_t$.
\end{theo}

Moreover, in some cases our technique allows us to strengthen the previous theorem and  to
obtain a prescribed {\em $1$-factor}, i.e., a perfect packing by cycles with given lengths.
To illustrate this we prove the following result, in which we also assume that $G$ is a
tournament to make the proof more convenient to present, although one can
remove this assumption.

\begin{theo} \label{1-factor}
For any number $M$ there is $c>0$ and numbers $T$ and $n_0$ so that
if $G$ is a tournament on $n>n_0$ vertices
with minimum indegree and outdegree at least $(1/2-c)n$
and $n_1,\cdots,n_t$ are numbers satisfying $\sum_{i=1}^t n_i = n$
then $G$ contains a $1$-factor with cycle lengths
$n_1,\cdots,n_t$ if the following holds:
for some $3 \le k \le M$ at least $T \log n$
of the $n_i$ are equal to $k$ and at least $T$ of the $n_i$
lie between $k+1$ and $M$.
\end{theo}

The rest of this paper is organised as follows. In the next
section we prove Theorem \ref{triangles} using probabilistic arguments
(R\"odl nibble) together with
the idea of `absorbing structures' introduced by R\"odl, Rucinski and Szemer\'edi.
In section 3 we prove Theorem \ref{long-cycles}, which is
the first ingredient in the proof of Theorem \ref{almost-factor},
and a result of independent interest: an asymptotically best
possible condition for finding a $1$-factor in which all prescribed
cycle lengths are long. To deal with short cycles we need
the machinery of Szemer\'edi's Regularity Lemma and the blowup lemma
of Koml\'os, S\'ark\"ozy and Szemer\'edi, which we describe in section 4.
Then in section 5 we prove Theorem \ref{ck}, the second ingredient in
the proof of Theorem \ref{almost-factor}, giving an almost perfect
packing by $k$-cycles for any fixed $k$. Section 6 contains the
proofs of Theorems \ref{almost-factor} and \ref{1-factor} and
the final section contains some concluding remarks.

\medskip

\nib{Notation.} Given two vertices $x$ and $y$ of a directed graph $G$, we
write $xy$ for the edge directed from $x$ to $y$.  We write $N^+_G(x)$
for the outneighbourhood of a vertex $x$ and $d^+_G(x):=|N^+_G(x)|$
for its outdegree.  Similarly, we write $N^-_G(x)$ for the
inneighbourhood of $x$ and $d^-_G(x):=|N^-_G(x)|$ for its indegree.
We write $N_G(x):=N^+_G(x)\cup N^-_G(x)$ for the neighbourhood of $x$
and $d_G(x):=|N_G(x)|$ for its degree. We use $N^+(x)$ etc. whenever this is
unambiguous. As is customary in Extremal Graph Theory, our approach to the
problems researched will be asymptotic in nature. We thus assume
that the order $n$ of a graph $G$ tends to infinity and therefore is sufficiently
large whenever necessary. We also assume that the constant $c$, which
controls the deviation of the degrees of $G$ from $n/2$ is
sufficiently small.
When we speak of `paths' and `cycles' in directed graphs
it is always to be understood that these are directed paths and
cycles. We use the notation $0 < \alpha \ll \beta$ to mean that
there is an increasing function $f(x)$ so that the following argument
is valid for $0 < \alpha < f(\beta)$.
We write $a \pm b$ to
denote an unspecified real number in the interval $[a-b,a+b]$.

\section{Covering by cyclic triangles}

In this section we prove Theorem \ref{triangles}. Our arguments combine
probabilistic reasoning (R\"odl nibble) together with idea of
 `absorbing structures' introduced by R\"odl, Rucinski and Szemer\'edi.
We divide the exposition into five subsections, that successively
treat two simple lemmas, large deviation inequalities,
the nibble, our absorbing structure, and the proof of Theorem \ref{triangles}.

\subsection{Two simple lemmas}

Our first lemma shows that, under the hypotheses of Theorem \ref{triangles},
there are approximately the same number of cyclic triangles through
every vertex.

\begin{lemma} \label{count-tri}
Suppose $c>0$ and $G$ is an oriented graph on $n$ vertices
with minimum indegree and outdegree at least $(1/2-c)n$.
Then every vertex $x$ of $G$ belongs to at least $\big(1/8 - 2c\big)n^2$ and at most
$\big(1/8 + 2c\big)n^2$ cyclic triangles.
\end{lemma}

\nib{Proof.}
To prove this lemma, we need to estimate $e(N^+(x),N^-(x))$, the number of edges
in $G$ going from the outneighbourhood of $x$ to the inneighbourhood
of $x$.
\begin{align*}
e(N^+(x),N^-(x)) & \geq \sum_{y \in N^+(x)} \Big (d^+(y) - |V(G) \sm (N^+(x) \cup N^-(x))|\Big) - e(N^+(x)) \\
& \geq |N^+(x)|\Big((1/2-c)n - n + |N^+(x)| + |N^-(x)|\Big) - |N^+(x)|^2/2 \\
& = |N^+(x)|\Big((1/2-c)n - n + |N^+(x)|/2 + |N^-(x)|\Big)\\
& \geq (1/2-c)n\Big((1/2-c)n - n + \frac{3}{2}(1/2-c)n\Big) \\
& = (1/2-c)n(1/4-5c/2)n \geq (1/8-2c)n^2.
\end{align*}
By symmetry we can also estimate $e(N^-(x),N^+(x)) \geq (1/8-2c)n^2$.
Therefore
$$\hspace{0.5cm}
e(N^+(x),N^-(x)) \leq |N^+(x)||N^-(x)| - e(N^-(x),N^+(x)) \leq n^2/4 - (1/8-2c)n^2 = (1/8+2c)n^2. \hspace{0.5cm} \Box $$

\medskip

Our next lemma will allow us to find cyclic triangles on `most' pairs
of vertices. Suppose that $G$ is an oriented graph on $n$ vertices.
We say that an edge $e$ of $G$ is $a$-good
if there are at least $a n$ cyclic triangles containing $e$;
otherwise we say it is $a$-bad. Also, given a vertex $x$ we
say that a vertex $y$ is $a$-good for $x$ if
an edge between $x$ and $y$ is $a$-good; otherwise we say it is $a$-bad for $x$.

\begin{lemma} \label{bad}
Suppose $c>0$ and $G$ is an oriented graph on $n$ vertices
with minimum indegree and outdegree at least $(1/2-c)n$.
For any $a>0$ and vertex $x$, there are at most
$(2a + 4c)n$ $a$-bad vertices for $x$ in each of
$N^+(x)$ and $N^-(x)$, so the total number of
$a$-bad vertices for $x$ is at most $(4a + 10c)n$.
\end{lemma}

\nib{Proof.}
Let $S$ be the set of $a$-bad vertices for $x$ that
belong to $N^+(x)$. Then, by definition, any $y \in S$ has
at most $a n$ outneighbours in $N^-(x)$. By averaging
there is some $y \in S$ with at most $|S|/2$ outneighbours
in $S$, and for this $y$ we have
\begin{align*}
(1/2-c)n & \leq |N^+(y)| = |N^+(y) \cap N^-(x)| + |N^+(y) \cap N^+(x)|
+ |N^+(y) \sm N(x)| \\
& \leq an+\big(|N^+(x)|-|S|/2\big)+|V(G) \sm N(x)| = an -|S|/2 +n - |N^-(x)|\\
& \leq an + (1/2+c)n - |S|/2,
\end{align*}
so $|S| \leq  (2a + 4c)n$. Similarly there are
at most $(2a + 4c)n$ $a$-bad vertices for $x$ that
belong to $N^-(x)$. Since $|V(G) \sm N(x)| \le 2cn$ there are at most
$(4a + 10c)n$ $a$-bad vertices for $x$. \qed

\subsection{Large deviation inequalities}

We will use the following three large deviation estimates. The first one is a
classical Chernoff-type bound see, e.g., \cite{AS} Appendix A.

\begin{theo} \label{chernoff}
Suppose $a>0$ and $X_1,\cdots,X_m$ are independent identically
distributed random variables with
$\mathbb{P}(X_i=1)=p$ and $\mathbb{P}(X_i=0)=1-p$. Then

1. $\mathbb{P}(\sum_{i=1}^m X_i < pm - a) < e^{-a^2/2pm}$, and

2. $\mathbb{P}(\sum_{i=1}^m X_i \ge pm + a) < e^{-a^2/2pm + a^3/2(pm)^2}$.

\end{theo}

\noindent
Another useful inequality is due to Azuma (see, e.g, \cite{AS, JLR}).
\begin{theo} \label{azuma}
Let $t_1, \ldots, t_n$ be a family of independent indicator random variables.
Suppose that real-valued function $X=X(t_1, \dots, t_n)$ is $c$-Lipschitz, i.e., changing the value of
any $t_i$ can change the value of $X$ by at most $c$. Then
$$\mathbb{P} \big( |X-\mathbb{E}X| >a\big) \leq 2e^{-a^2/2nc^2}\,.$$
\end{theo}

The third inequality we need was proved by Kim and Vu
(see, e.g., Chapter 7.8 in \cite{AS}).
Suppose $f(x_1,\cdots,x_m) = \sum_{e \in H} \prod_{i \in e} x_i$
is a homogeneous polynomial of degree $k$ defined by a $k$-uniform hypergraph $H$
on $[m]$. Let $X_1,\cdots,X_m$ be independent identically
distributed random variables with
$\mathbb{P}(X_i=1)=p$, $\mathbb{P}(X_i=0)=1-p$ and
let $Y = f(X_1,\cdots,X_m)$. For $J \subset [m]$ define
$\partial_J f$ to be the partial derivative of $f$ with respect to
the variables $\{x_i\}_{i \in A}$ and
the {\em non-zero influence} $\mb{E}'Y = \max_{|A| \ge 1} \mb{E} \partial_A f$.

\begin{theo} \label{kim-vu}
For $t>1$,
$$\mb{P}(|Y-\mb{E}Y| > (2k)!t^k\sqrt{\mb{E}Y\mb{E}'Y})
\le 16e^{-t+(k-1)\log m}.$$
\end{theo}

\subsection{The nibble}

The `nibble' is a term referring to a semi-random construction method,
used by R\"odl \cite{R} in proving the existence of asymptotically good designs.
Several researchers realised that this method applies in a more general setting,
dealing with matchings in uniform hypergraphs.
The following theorem is due to Pippenger, following Frankl and R\"odl,
with further refinements by Pippenger and Spencer \cite{PS}.
We refer the reader to the presentation given in \cite{F} and in
\cite{AS}, pp. 54--58. For a pair of vertices $x, y$ of a hypergraph
$H$, the common degree $d(x,y)$ is the number of edges of $H$ containing both $x$ and $y$.

\begin{theo} \label{pip}
For any $\eps>0$ and number $r$ there is $\delta>0$ and a
number $d$ so that the following holds for $n>D>d$.
Any $r$-uniform hypergraph $H$ on $n$ vertices such that
any vertex $x$ has degree $d(x) = (1 \pm \delta)D$ and
any pair of vertices $x,y$ has common degree $d(x,y) < \delta D$
contains a matching covering at least $(1-\eps)n$ vertices.
\end{theo}

The result which appears in \cite{AS} deals with covering of vertices of hypergraph
rather than matchings.
It states that a $k$-uniform hypergraph $H$ on $n$ vertices with all degrees $(1 + o(1))D$
and all codegrees $o(D)$ has a collection of $(1+o(1))n/k$ edges which covers all its vertices.
It is easy to check that deletion of all pairs of intersecting edges from this collection
gives a matching covering $(1-o(1))n$ vertices of $H$.

We also need a further property of the matching
which comes out of the proof of Theorem \ref{pip} by means
of a semi-random `nibble'. The matching is constructed
in a series of `bites', in which we choose each remaining available
edge independently with some probability $p=\Theta(n^{-2})$
(which shrinks by a constant factor with each step)
and delete any pair of edges that intersect. It is shown
that with probability at least $0.9$ (say)
each bite preserves certain regularity
properties in the hypergraph that allow the nibble to proceed.
The parameters of the proof are such that we
may assume that the first bite constructs a matching
of size $\beta n$ with $\eps \ll \beta \ll 1$.

\subsection{The absorbing structure}

The nibble can be used to cover all but $o(n)$ vertices, and to
make further progress we will need a mechanism that will allow us
to gradually `absorb' the remaining vertices into our triangle packing.
Our approach was inspired by ideas used in \cite{RRS1, RRS2} to obtain
results on matchings and Hamiltonian hypergraphs.

\parbox{0.55\textwidth}{
Suppose $G$ is an oriented graph and $Q = \{v_1,v_2,v_3,v_4\}$
is a quadruple of vertices in $G$. We say that the disjoint sets
$a_1 a_2 a_3$, $b_1 b_2 b_3$, $c_1 c_2 c_3$ are an
absorbing triple of triangles for $Q$ if each of the
following triples is a cyclic triangle in $G$:
$a_1 a_2 a_3$, $b_1 b_2 b_3$, $c_1 c_2 c_3$,
$v_1 a_1 b_1$, $v_2 c_1 a_2$, $v_3 b_2 c_2$,
$v_4 a_3 b_3$.
The motivation for this definition
is that if we have a set of disjoint cyclic triangles
$C_1, \cdots, C_t$ that includes
$a_1 a_2 a_3$, $b_1 b_2 b_3$, $c_1 c_2 c_3$
and is disjoint from $Q$, then we can enlarge our
collection by replacing $a_1 a_2 a_3$, $b_1 b_2 b_3$, $c_1 c_2 c_3$
by $v_1 a_1 b_1$, $v_2 c_1 a_2$, $v_3 b_2 c_2$,
$v_4 a_3 b_3$. Thus $Q$ is absorbed and the vertex
$c_3$ is lost, for a net gain of one triangle.
}\hspace{0.4cm}
\parbox{0.37\textwidth}{
\includegraphics[height=7cm]{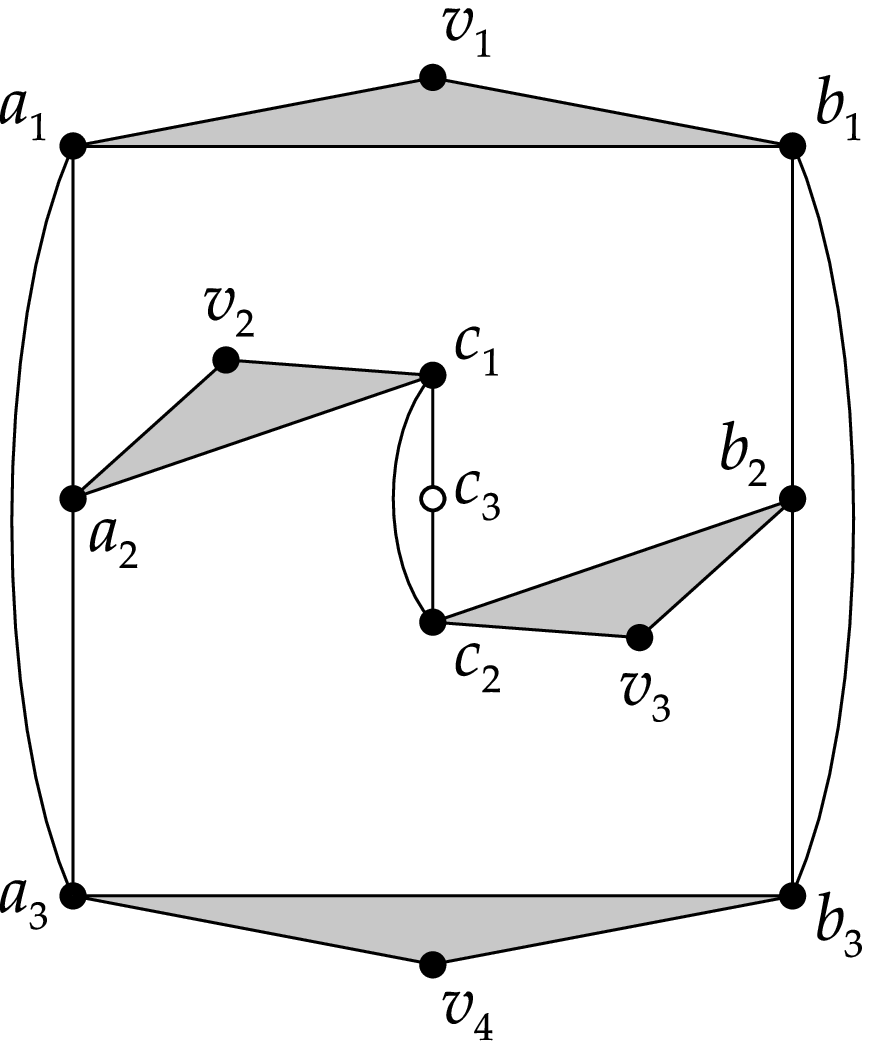}
}

The following lemma shows that there are many absorbing triples
for every quadruple $Q$.

\begin{lemma} \label{pre-absorb-tri}
There is some $c>0$ and number $n_0$ such that
if $G$ is an oriented graph on $n>n_0$ vertices
with minimum indegree and outdegree at least $(1/2-c)n$
then for any quadruple of vertices $Q$ there are at
least $(n/100)^9$ absorbing triples for $Q$ in $G$.
\end{lemma}

\nib{Proof.}
We use the above notation and greedily construct
the absorbing triples by repeated application of Lemma \ref{bad}.
\begin{enumerate}
\item
Pick $a_1$ to be $1/8$-good for $v_1$ and disjoint from $Q$.
There are at least $n - (4/8+10c)n - 4 > (1/2-11c)n$
possible choices.
\item
Pick $a_2$ to be $1/16$-good for $a_1$ and $v_2$
(and disjoint from $Q \cup \{a_1\}$: we will not keep
repeating this condition). There are at least
$n - 2(4/16+10c)n -5 > (1/2 - 21c)n$ possible choices.
\item
Pick $a_3$ to be $1/128$-good for $v_4$ and
so that $a_1 a_2 a_3$ is a cyclic triangle. Since
$a_2$ is $1/16$-good for $a_1$ there are at least
$n/16 - (4/128+10c)n - 6 > (1/32-11c)n$ possible choices.
\item
Pick $b_3$ so that $v_4 a_3 b_3$ is a cyclic triangle.
Since $a_3$ is $1/128$-good for $v_4$
there are at least $n/128-7$ possible choices.
\item
Pick $b_1$ to be $1/64$-good for $b_3$ and so that
$v_1 a_1 b_1$ is a cyclic triangle. Since
$a_1$ is $1/8$-good for $v_1$
there are at least $n/8 - (4/64+10c)n-8 > (1/16-11c)n$ possible choices.
\item
Pick $b_2$ to be $1/512$-good for $v_3$ and so that
$b_1 b_2 b_3$ is a cyclic triangle. Since $b_1$ is $1/64$-good for $b_3$
there are at least $n/64 - (4/512+10c)n-9 > (1/128-11c)n$ possible choices.
\item
Pick $c_2$ so that $v_3 b_2 c_2$ is a cyclic triangle. Since
$b_2$ is $1/512$-good for $v_3$
there are at least $n/512-10$ possible choices.
\item
Pick $c_1$ to be $1/128$-good for $c_2$ and so that
$v_2 c_1 a_2$ is a cyclic triangle. Since
$a_2$ is $1/16$-good for $v_2$
there are at least $n/16 - (4/128+10c)n-11 > (1/32-11c)n$ possible choices.
\item
Pick $c_3$ so that $c_1 c_2 c_3$ is a cyclic triangle. Since
$c_1$ is $1/128$-good for $c_2$
there are at least $n/128-12$ possible choices.
\end{enumerate}
Note that given three cyclic triangles there are $3!=6$ ways to chose which
one of them is going to be $a_1 a_2 a_3$, $b_1 b_2 b_3$, $c_1 c_2 c_3$.
Then for every cyclic triangle there are $3!=6$ different labeling
of its vertices. This implies that each configuration of three cyclic triangles
was counted at most $6^4$ times. Hence the number of absorbing triples is at least
\begin{eqnarray*}
&~&6^{-4} \cdot \Big(\frac{1}{2}-11c\Big)n \cdot \Big(\frac{1}{2} - 21c\Big)n \cdot
\Big(\frac{1}{32}-11c\Big)n
\cdot \Big(\frac{n}{128}-7\Big) \cdot \Big(\frac{1}{16}-11c\Big)n \cdot
\Big(\frac{1}{128}-11c\Big)n\\
&~&\cdot \Big(\frac{n}{512}-10\Big)
\cdot \Big(\frac{1}{32}-11c\Big)n \cdot \Big(\frac{n}{128}-12\Big)
= \big(6^{-4}2^{-46}-O(c)\big)n^9 > \big(n/100\big)^9.  \hspace{2.7cm} \Box
\end{eqnarray*}

\medskip

Next we use the previous lemma to show that a random selection of vertex-disjoint
cyclic triangles will have many absorbing triples for every quadruple of vertices of
$G$.

\begin{lemma} \label{absorb-tri}
Suppose $0<1/n_0 \ll c \ll c_2 \ll 1$ and
$G$ is an oriented graph on $n>n_0$ vertices
with minimum indegree and outdegree at least $(1/2-c)n$.
Suppose we form a collection of vertex-disjoint
cyclic triangles $C$ by choosing each cyclic triangle in $G$
independently with probability $p=c_2/n^2$
and deleting any pair of triangles that intersect.
Then with probability at least $0.9$ we have
$|C| = m = (1 \pm c_2^{1/2})c_2n/24$ and for any quadruple of vertices $Q$
there are at least $10^{-16} m^3$ absorbing triples for $Q$ in $C$.
\end{lemma}

\nib{Proof.}
Let $C'$ be a collection of cyclic triangles formed by choosing each cyclic triangle of $G$
randomly and independently with probability $p=c_2/n^2$.
By Lemma \ref{count-tri}, every vertex of $G$ is contained in
$(1/8\pm 2c)n^2$ cyclic triangles and therefore the number of cyclic
triangles in $G$ is
$T=(1/8\pm 2c)n^2 \cdot n/3= (1 \pm 16c)n^3/24$.
Applying Chernoff bounds (mentioned in Section 2.2) we obtain that
$$|C'|= (1 \pm c)p^3T= (1 \pm 20c)pn^3/24 = (1 \pm 20c)c_2n/24$$ with high probability.
Let $Z$ be the number of pairs of intersecting triangles in $C'$.
Since the total number of such pairs is clearly at most $n^5$ we have that
$\mb{E}Z < p^2 n^5 = c_2^2n$. Hence, $c_2 \ll 1$ together with Markov's inequality
gives that $Z < c_2^{3/2} n/100$ with probability at least
$0.95$. Since $c \ll c_2$,
by definition of $C$, we obtain
$$m=|C| \ge |C'|-2Z > (1 \pm c_2^{1/2})c_2n/24.$$

Given a quadruple of vertices $Q$, let $A_Q$ be the set
of absorbing triples for $Q$. By Lemma \ref{pre-absorb-tri} we have
$|A_Q| > (n/100)^9$. Let $X_Q$ be the random variable counting the number of absorbing
triples for $Q$ that belong to $C'$. Then
$\mb{E} X_Q = p^3|A_Q| > (c_2 n)^3/10^{18}$. We can write
$X_Q = \sum_{S \in A_Q} \prod_{T \in S} I_T$
where $I_T$ is the indicator random variable for the event
that triangle $T$ is chosen for $C'$.
Since $X_Q$ is a homogeneous polynomial of degree three,
we can estimate the probability that $X_Q$ is small
by Theorem \ref{kim-vu}.
Since the number of absorbing triples for $Q$ which contains a given triangle is clearly at most $n^6$
it is easy to see that $\max_{|J|= 1} \mb{E} \partial_J X_Q$ is at most
$p^2 n^6 = (c_2 n)^2$. Similarly, there are at most $n^3$
absorbing triples containing a given pair of triangles and therefore
$\max_{|J|= 2} \mb{E} \partial_J
X_Q \leq p n^3=c_2 n$. This implies that the non-zero influence
$\mb{E}' X_Q = \max_{|J| \ge 1} \mb{E} \partial_J X_Q$
is bounded by $p^2 n^6 = (c_2 n)^2$. Thus, choosing
$t = \Theta(pn^3)^{1/6} = \Theta(n^{1/6})$
in Theorem \ref{kim-vu}
we can estimate
$$\mb{P}(X_Q < \mb{E}X_Q/2) \le 16e^{-t+2\log n^3} \ll n^{-4}.$$
Taking a union bound over all quadruples of vertices of $G$, we obtain that
with high probability $X_Q > p^3|A_Q|/2 > m^3/10^{15}$ for every $Q$.
Note also that deletion of any triangle from $C'$ can destroy at most
$|C'|^2 <2 m^2$ absorbing triples for $Q$. Since we delete at most
$2Z < c_2^{3/2}n/50 < c_2 m/2$
triangles to form $C$  and since $c_2 \ll 1$, we still have at least
$ m^3/10^{15}- c_2 m^3> m^3/10^{16}$ absorbing triples for
each $Q$.
\qed

\subsection{Proof of Theorem \ref{triangles}}

(i) Choose constants to satisfy the hierarchy
$0 < 1/n \ll c \ll c_1 \ll c_2 \ll 1$
and suppose $G$ is an oriented graph on $n$ vertices
with minimum indegree and outdegree at least $(1/2-c)n$.
Consider the hypergraph $H$ on the same vertex set of $G$ whose
edges are all cyclic triangles in $G$. By Lemma \ref{count-tri}
every vertex $x$ in $H$ has degree $d_H(x) = (1/8 \pm 2c)n^2$.
Also, for any pair of vertices $x,y$ we have $d_H(x,y) \le n-2 \ll n^2/8$.
Applying Theorem \ref{pip}, we can cover all but at most $(1-c_1)n$ vertices with
vertex-disjoint cyclic triangles.

Furthermore, as explained in the paragraph after Theorem \ref{pip},
we may assume that the first bite of the nibble was obtained by
choosing each cyclic triangle in $G$ with probability $c_2/n^2$ and deleting any pair of triangles
that intersect.
Since the bite was valid for the nibble with probability at least $0.9$, we can
also assume that it is an absorbing collection $C$ as given by Lemma \ref{absorb-tri}.
Now, as long as there at least $4$ uncovered vertices, we can repeatedly choose
a quadruple $Q$ from these vertices and increase our triangle packing by
using an absorbing triple for $Q$ from $C$. Note that at each such iteration,
we can only use absorbing triples from $C$ no triangle of which has yet been used in previous rounds.
Since at each iteration we use three triangles and each triangle can participate in at
most $|C|^2=m^2$ absorbing triples for $Q$, we destroy at most
$3m^2$ absorbing triples for $Q$ at every round.
As the number of rounds is at most $c_1n \ll m = \Omega(c_2n)$ there are still at least
$10^{-16}m^3-c_1n m^2>10^{-17}m^3$ absorbing triples remaining untouched during the whole procedure.
This shows that our process can be continued until only 3 vertices will remain uncovered,
which completes the proof of the first part.

\noindent (ii) To prove the second part of the theorem,
partition a set of $n$ vertices with $n=18k+3$
into three sets $V_0,V_1,V_2$ of size
$|V_0|=6k$, $|V_1|=6k+1$, $|V_2|=6k+2$.
Construct a tournament $T$ as follows. Between the classes we orient all pairs
from $V_i$ to $V_{i+1}$, where addition is mod $3$. Inside each
class we place a tournament that is as regular as possible, i.e.,
in $V_1$ all indegrees and outdegrees are $3k$, in $V_0$ each vertex has
indegree and outdegree $3k$ and $3k-1$ in some order, and in $V_2$
each vertex has indegree and outdegree $3k$ and $3k+1$ in some order.
Then every vertex in $T$ either has indegree and outdegree $9k+1$
or indegree and outdegree $9k$ and $9k+2$ in some order. However,
any  collection of vertex-disjoint cyclic triangles in $T$ must
leave at least $3$ vertices uncovered. To see this, notice
that a cyclic triangle must either have one point in each part
or all three points in one of $V_i$, and so however many triangles
we remove from $T$ the class sizes will always be different mod $3$. \qed

\section{Long cycles}

Our first ingredient in the proof of Theorem \ref{almost-factor} will be the
following theorem, which shows that the minimum semidegree threshold
for finding a $1$-factor in which all the prescribed cycle lengths are
large is asymptotically $3n/8$. The lower bound is given by a
construction in \cite{H} (see also \cite{KelKO}) of an oriented graph with minimum semi-degree
$\sim 3n/8$ and with no $1$-factor at all. Hence it only remains to prove the upper bound.

\begin{theo} \label{long-cycles}
For any $\delta>0$ there are numbers $M$ and $n_0$ so that
if $G$ is an oriented graph on $n>n_0$ vertices
with minimum indegree and outdegree at least $(3/8+\delta)n$
and $n_1,\cdots,n_t$ are numbers satisfying
$n_i \ge M$ for $1 \le i \le t$ and $\sum_{i=1}^t n_i = n$
then $G$ contains a $1$-factor with cycle lengths
$n_1,\cdots,n_t$.
\end{theo}

Our proof combines the partitioning argument
similar to that used in \cite{AF} together with the result of \cite{KelKO2}.
In Theorem 8 of \cite{KelKO2} it was proved that for all sufficiently large
$n$ every oriented graph $G$ on $n$
vertices with minimum semidegree at least $(3n-4)/8$ contains an $ \ell$-cycle for
all $3 \leq \ell \leq n$.
We also need another large deviation inequality, for
the hypergeometric random variable $X$ with parameters $(n,m,k)$, which is
defined as follows. Fix $S \subset [n]$ of size $|S|=m$. Pick
a random $T \subset [n]$ of size $|T|=k$. Define $X=|T \cap S|$.
Then $\mb{E}X = km/n$. We have the following
`Chernoff bound' approximation for $0<a<3/2$
(see \cite{JLR} pp. 27--29):
\begin{equation} \label{hyp-geom}
\mb{P}(|X - \mb{E}X| > a\mb{E}X) < 2 e^{-\frac{a^2}{3}\mb{E}X}.
\end{equation}

The following lemma is an immediate consequence of the previous inequality.

\begin{lemma} \label{partition}
For any $\alpha, \beta >0$ there is a number $n_0$ so that the following holds.
Suppose $G$ is an oriented graph on $n>n_0$ vertices
with minimum indegree and outdegree at least $\alpha n$
and $\beta n < m < (1-\beta)n$. Then there is a partition of $V(G)$ as
$A \cup B$ with $|A|=m$ and $|B|=n-m$ so that
$G[A]$ has minimum indegree and outdegree at least $(\alpha-n^{-1/3})m$ and
$G[B]$ has minimum indegree and outdegree at least $(\alpha-n^{-1/3})(n-m)$.
\end{lemma}

\nib{Proof of Theorem \ref{long-cycles}.}
Without loss of generality suppose that $n_1 \ge n_i$ for $1 \le i \le t$.
Consider two cases.

\nib{Case 1.}
If $n_1 > (1-\delta/2)n$ then $\sum_{i\geq 2} n_i < \delta n/2$ and we can use Theorem 8 in
\cite{KelKO2} (mentioned above)
to choose disjoint cycles of length $n_2,\cdots,n_t$ one by one. Indeed, it is possible since
during this process the semidegree of the oriented graph which remains is always
at least $(3/8 +\delta)n -\sum_{i\geq 2} n_i > 3n/8$. In particular,
the oriented graph on $n_1$ vertices which we obtain in the end has minimum semidegree
larger than $3n/8 \geq 3n_1/8$ and therefore has a Hamilton cycle, so we are done.
Note that this argument works as long as the minimum semidegree of the graph is at least
$(3/8 +\delta/2)n$, which will be used in the analysis of the second case.

\nib{Case 2.} If $n_1 \leq (1-\delta/2)n$ we can partition
$[t]=I \cup J$ so that $n_I = \sum_{i \in I} n_i$ and
$n_J = \sum_{i \in J} n_i$ are both at most $(1-\delta/2)n$
(we may assume $\delta<1/3$). Then by Lemma \ref{partition} there is a
partition of the vertices of $G$ into sets $V_I$ of size $n_I$ and
$V_J$ of size $n_J$ so that
$G[V_I]$ has minimum semidegree at least $(3/8+\delta-n^{-1/3})n_I$ and
$G[V_J]$ has minimum semidegree at least $((3/8+\delta-n^{-1/3})n_J$.

Now we repeat the above splitting procedure for both $V_I$ and $V_J$, repeatedly
partitioning while Case 2 holds. Each time the number of vertices in the part which
was split is reduced by a factor of $(1-\delta/2)$ and no part in our process ever has
size smaller than $M$. Therefore, for any part $S$ in the final partition the induced graph
$G[S]$ has minimum indegree and outdegree at least
$$\left( 3/8+\delta - M^{-1/3} \sum_{i=0}^\infty (1-\delta/2)^{i/3} \right)
|S|
= \left( 3/8+\delta -  \frac{1}{1-(1-\delta/2)^{1/3}} M^{-1/3} \right) |S|.$$
For large enough $M$ this is more than $(3/8+\delta/2)|S|$, so we can use the argument of Case 1
to find the required cycles. \qed

\section{Regularity}

The second ingredient in the proof of Theorem \ref{almost-factor} will be
a theorem giving an almost perfect packing of $k$-cycles when $k$ is fixed
and the number of vertices $n$ is large. The proof of this theorem will use
the machinery of Szemer\'edi's Regularity Lemma and the blowup lemma
of Koml\'os, S\'ark\"ozy and Szemer\'edi, which we will now describe. We
will be quite brief, so for more details and motivation we refer the reader
to the surveys \cite{KSi} for the regularity lemma and \cite{Ko} for the blowup lemma.

We start with some definitions. The density of a bipartite graph $G = (A,B)$
with vertex classes $A$ and $B$ is defined to be
\[d_G(A,B) := \frac{e_G(A,B)}{|A||B|}.\]
We often write $d(A,B)$ if this is unambiguous. Given $\eps>0$, we say
that $G$ is \emph{$\eps$-regular} if
for all subsets $X\subseteq A$ and $Y\subseteq B$
with $|X|>\eps|A|$ and $|Y|>\eps|B|$ we have that
$|d(X,Y)-d(A,B)|<\eps$. Given $d\in [0,1]$ we say
that $G$ is $(\eps,d)$-\emph{super-regular}
if it is $\eps$-regular and furthermore $d_G(a)\ge (d-\eps) |B|$
for all $a\in A$ and
$d_G(b)\ge (d-\eps)|A|$ for all $b\in B$.
(This is a slight variation of the standard definition
of $(\eps,d)$-super-regularity where one requires $d_G(a)\ge d |B|$ and
$d_G(b)\ge d|A|$.)

The Diregularity Lemma is a version of the Regularity Lemma for digraphs due to
Alon and Shapira \cite{ASh} (with a similar proof to the undirected version).
We will use the following degree form of the Diregularity Lemma, which can be easily derived
(see e.g. \cite{Yo}) from the standard version, in exactly the same manner
as the undirected degree form.

\begin{lemma}[Degree form of the Diregularity Lemma]\label{direg}
For every $\eps \in (0,1)$ and $M'>0$ there are numbers $M$ and $n_0$
such that if $G$ is a digraph on $n\ge n_0$ vertices and
$d\in[0,1]$,
then there is a partition of the vertices of $G$ into
$V_0,V_1,\cdots,V_s$ and a spanning subdigraph $G'$ of $G$ such that
 the following holds:
\begin{itemize}
\item $M'\le s\leq M$,
\item $|V_0|\leq \eps n$,
\item $|V_1|=\cdots=|V_s|$,
\item $d^+_{G'}(x)>d^+_G(x)-(d+\eps)n$ for all vertices $x\in G$,
\item $d^-_{G'}(x)>d^-_G(x)-(d+\eps)n$ for all vertices $x\in G$,
\item for all $i=1,\cdots,s$ the digraph $G'[V_i]$ is empty,	
\item for all $1\leq i,j\leq s$ with $i\neq j$
the bipartite graph whose vertex classes
are $V_i$ and $V_j$ and whose edges are all the edges in $G'$ directed from $V_i$ to $V_j$
is $\eps$-regular and has density either $0$ or density at least $d$.
\end{itemize}
\end{lemma}

Given clusters $V_1,\cdots,V_s$ and a digraph $G'$,
the \emph{reduced digraph $R'$
with parameters $(\eps,d)$} is the digraph
whose vertex set is $[s]$ and in which $ij$ is an edge
if and only if the bipartite graph whose vertex classes
are $V_i$ and $V_j$ and whose edges are all the edges in $G'$ directed from $V_i$ to $V_j$
is $\eps$-regular and has density at least $d$.
(So $ij$ is an edge in $R'$ if and only if there is an edge
from $V_i$ to $V_j$ in $G'$.)
It is easy to see that the reduced digraph $R'$ obtained
from the regularity lemma
`inherits' the minimum degree of $G$, in that
$\delta^+(R')/|R'| > \delta^+(G)/|G| - d - 2\eps$
and $\delta^-(R')/|R'| > \delta^-(G)/|G| - d - 2\eps$.
However, $R'$ is not necessarily oriented even if the original digraph $G$ is.
The next lemma from \cite{KelKO}
shows that by discarding edges with appropriate
probabilities one can go over to a
reduced oriented graph $R\subseteq R'$ which still inherits the minimum degree
and density of $G$.

\begin{lemma} \label{red-orient}
For every $\eps\in (0,1)$ there exist numbers
$M'=M'(\eps)$ and $n_0=n_0(\eps)$
such that the following holds. Let $d\in [0,1]$ with $\eps \ll d$,
let $G$ be an oriented graph of order
$n \ge n_0$ and let $R'$ be the reduced digraph with parameters $(\eps,d)$
obtained by applying the Lemma \ref{direg} to $G$ with parameters
$\eps, d$ and $M'$.
Then $R'$ has a spanning oriented subgraph $R$ such that
$\delta^+(R)\ge (\delta^+(G)/|G|-(d+3\eps))|R|$ and
$\delta^-(R)\ge(\delta^-(G)/|G|-(d+3\eps))|R|$.
\end{lemma}

We conclude this section with the blow-up lemma of Koml\'os, S\'ark\"ozy and Szemer\'edi \cite{KSS1}.

\begin{lemma} \label{blowup}
Given a graph $F$ on $[s]$ and positive numbers $d,\Delta$, there is
a positive real $\eta_0=\eta_0(d,\Delta,s)$
such that the following holds for all positive numbers
$\ell_1,\dots,\ell_s$ and all $0<\eta\le \eta_0$. Let $F'$ be the graph obtained
from $F$ by replacing each vertex $i\in F$ with a set $V_i$ of $\ell_i$ new
vertices and joining
all vertices in $V_i$ to all vertices in $V_j$ whenever $ij$ is an edge
of $F$. Let $G'$ be a spanning subgraph of $F'$
such that for every edge $ij\in F$ the bipartite graph consisting of
all the edges of $G'$ between the sets $V_i, V_j$ is
$(\eta,d)$-super-regular. Then $G'$ contains a copy of every subgraph $H$ of
$F'$ with $\Delta(H)\le \Delta$. Moreover, this copy of $H$ in $G'$
maps the vertices of $H$
to the same sets $V_i$ as the copy of $H$ in $F'$, i.e. if $h \in V(H)$
is mapped to $V_i$ by the copy
of $H$ in $F'$, then it is also mapped to $V_i$ by the copy of $H$ in $G'$.
\end{lemma}

Note that the `moreover' part of this statement does not appear in the usual formulation
of the Blow-up Lemma but is stated explicitly in its proof.

\section{Short cycles}

We now come to the second ingredient in the proof of Theorem \ref{almost-factor}, which is
the following statement, providing an almost perfect packing by $k$-cycles,
when $k$ is fixed and $n$ is large.

\begin{theo} \label{ck}
For any number $k \ge 3$ there is some real $c>0$ and numbers $C$ and $n_0$
so that if $G$ is an oriented graph on $n>n_0$ vertices
with minimum indegree and outdegree at least $(1/2-c)n$
then $G$ contains vertex-disjoint $k$-cycles covering all
but at most $C$ vertices.
\end{theo}

Since the case $k=3$ was already proved in Theorem \ref{triangles},
in the rest of this section we assume that
$k \geq 4$. First we note that the following result is an immediate
consequence of Lemma \ref{blowup} (the blow-up lemma).

\begin{coro} \label{tri-blowup}
Suppose $H$ is an oriented graph with parts $V_0,V_1,V_2$ of equal size, so that
all edges go from $V_i$ to $V_{i+1}$ (addition mod $3$).
Suppose also that the underlying graphs between each pair of classes are
$(\eta,d)$-super-regular, for some $\eta \ll d < 1$. Then $H$ has
a perfect packing by cyclic triangles (with one vertex in each class).
\end{coro}

Next we use this corollary to obtain an almost perfect packing by $k$-cycles
under an additional semidegree assumption on the parts.

\begin{theo} \label{ck-blowup}
Suppose numbers $k,m$ and reals $\delta,d,\eta$ satisfy
$0 < 1/m \ll \delta \ll \eta \ll d \ll 1/k \le 1/4$
and $H$ is an oriented graph whose vertices are partitioned into three parts $V_0,V_1,V_2$ of sizes
$0.9m \le |V_i| \le m$ satisfying

(i) $H[V_i]$ has minimum degree at least $(1-\delta)|V_i|$ for $i=0,1,2$,

(ii) the edges of $H$ between $V_i$ and $V_{i+1}$ are all directed from $V_i$ to $V_{i+1}$.

(iii) the underlying graphs between each pair of classes $(V_i,V_{i+1})$ are
$(\eta ,d)$-super-regular.

\noindent
Then $H$ contains a packing of $k$-cycles covering all but at most $3k$ vertices.
\end{theo}

\nib{Proof.}
First it will be useful
to see how to find such $k$-cycles in the oriented graph $K$ which is obtained from
$H$  by adding all directed edges from $V_i$ to $V_{i+1}$. We will
choose our $k$-cycles to have $k-2$ points in one class and $1$ point
in each of the other two classes. Let $n_i$ be the number of cycles
with $k-2$ points in $V_i$.  We need to choose $n_i$ so
that $|V(H)|/k-3 \le n_0 + n_1 + n_2 \le |V(H)|/k$ subject to the conditions
$|V_i| \ge (k-2)n_i + \sum_{j \ne i} n_j = (k-3)n_i + n_0+n_1+n_2$.
We may take $n_i = \lfloor (|V_i| - |V(H)|/k)/(k-3) \rfloor$.
Indeed, then $n_i= (|V_i| - |V(H)|/k)/(k-3)-x_i$ for some $0 \le x_i<1$.
With $ x=x_0+x_1+x_2$, we have
$$n_0+n_1+n_2=(|V_0|+|V_1|+|V_2|-3|V(H)|/k)/(k-3)-x_0-x_1-x_2 =
|V(H)|/k- x \,.$$
Since, by definition $0 \leq x <3$,
this implies that $|V(H)|/k-3 \le n_0 + n_1 + n_2 \le |V(H)|/k$,
and that
$$ |V_i|- (k-3)n_i = |V(H)|/k+(k-3)x_i=n_0+n_1+n_2+x+(k-3)x_i \geq n_0+n_1+n_2.$$
Since $k \geq 4$, we also have  $n_i > (0.9 - 3/k)m/(k-3) > m/4k$.
In order to form the cycles in oriented graph $K$ it clearly suffices to find $n_i$ disjoint directed
paths of length $k-2$ in $V_i$. Let $P_i$ be an arbitrary subset of $V_i$ of size
precisely $(k-2)n_i$. Ignore the direction of the edges and
consider the induced subgraph $H[P_i]$. By assumption (i) of the theorem,
the degree of every vertex in this graph is at least
$|P_i|-\delta|V_i|$. Since $|P_i|=(k-2)n_i=\Theta(m)$ and $\delta \ll 1/k$,
every vertex in $H[P_i]$ has at least
$(1-\frac{1}{k-2})|P_i|$ neighbours. Thus, applying the Hajnal-Szemer\'edi theorem
(mentioned in the introduction) we can find a collection of disjoint cliques of size $k-2$
covering all vertices of $P_i$. With the directions each of these cliques
becomes a tournament (i.e. a complete oriented graph) and it is well known
that any tournament has a Hamiltonian path (i.e. a path containing all of its vertices).
Given these paths, we may assign each path an arbitrary pair of vertices in $V_{i-1}$
and $V_{i+1}$ to form the required $k$-cycles in the oriented graph $K$.

Now we will use the blowup lemma to show that the same strategy
works even when the edges from $V_i$ to $V_{i+1}$ no longer form a
complete bipartite graph, but do form a super-regular pair. We start
by picking randomly disjoint sets $S_{i,j} \subset V_i$ of size
$|S_{i,j}| = n_j$ for all $ j \ne i \in \{0,1,2\}$.
For sufficiently large $m$,
by the Chernoff bound for hypergeometric distributions (see inequality (\ref{hyp-geom})), we can assume
that $|N^+(v) \cap S_{i+1,i}| > (d-\eta)n_i - m^{2/3}$ and
$|N^-(v) \cap S_{i-1,i}| > (d-\eta)n_i - m^{2/3}$ for every $i=0,1,2$ and $v \in V_i$.
Next let $P_i$ be an arbitrary subset of $V_i \sm \cup_{j \ne i} S_{i,j}$
of size $(k-2)n_i$. Using the same argument as in the previous paragraph we can find
disjoint paths $P_{i,1}, \cdots, P_{i,n_i}$ each of length $k-2$ covering all vertices of $P_i$.
Denote the first and last vertices of $P_{i,j}$ by $x_{i,j}$ and $y_{i,j}$,
respectively and let $X_i = \{x_{i,1}, \cdots, x_{i,n_i}\}$, $Y_i = \{y_{i,1}, \cdots, y_{i,n_i}\}$.
Note that these sets have size linear in $m$.

By regularity, there are at most $\eta |V_i|$ vertices $v \in V_i$ with
$|N^+(v) \cap X_{i+1}| < (d-2\eta)|X_{i+1}|$ and
at most $\eta |V_i|$ vertices $v \in V_i$ with
$|N^-(v) \cap Y_{i-1}| < (d-2\eta)|Y_{i-1}|$. Call such vertices {\em bad}.
For $i \ne j \in \{0,1,2\}$ define new sets $S'_{i,j}$
and paths $P'_{i,t}$ as follows.
Choose sets $B_{i,j} \sub S_{i,j}$, $j \ne i$ containing all bad vertices in $S_{i,j}$
with $|B_{i,j}|$ equal to the first number larger than $2\eta|V_i|$ that
is divisible by $k-2$. For every $i$ and $j \ne i$
choose a collection of paths $\{P_{i,\ell}\}, \ell \in C_{i,j}$
containing only good vertices such that the sets of indices
$C_{i,j}, j \ne i$ are disjoint and have size $|C_{i,j}|= |B_{i,j}|/(k-2)$.
Note that this is possible since the number of bad vertices is
at most $2\eta m$, the number of paths $P_{i,\ell}$ is $n_i=\Theta (m/k)$ and  $\eta \ll
1/k$. For $j \ne i$
remove the paths $\{P_{i,\ell}\}, \ell \in C_{i,j}$ from $P_i$, adding their vertices
to $S_{i,j}$, and replace the vertices lost from $P_i$ with $\cup_{j \ne i} B^j_i$.
Delete the vertices of $ B_{i,j}$ from $S_{i,j}$ and call the new set $S'_{i,j}$.
Note that the size of $S'_{i,j}$ is still $n_i$ and it now contains only good vertices.
Since there are at least $4\eta|V_i|$ vertices in $\cup_{j \ne i} B_{i,j}$ and $\delta \ll \eta$ we can
again use the same argument as above to find disjoint paths of length $k-2$ covering all vertices in
$\cup_{j \ne i} B_{i,j}$. Add these new paths instead of the paths $P_{i,\ell}, \ell \in \cup_{j \ne i}
C_{i,j}$ which were removed and call
the new collection of paths $P'_{i,1}, \cdots, P'_{i,n_i}$. Also, let
$X'_i = \{x'_{i,1}, \cdots, x'_{i,n_i}\}$
and $Y'_i = \{y'_{i,1}, \cdots, y'_{i,n_i}\}$
be the sets of first  and last vertices for the new collection of paths $\{P'_{i,t}\}$.

Now consider three new $3$-partite oriented graphs $H_0, H_1, H_2$ defined
as follows. The parts of $H_i$ are $S'_{i-1,i}$, $S'_{i+1,i}$ and an auxiliary
set of size $n_i$, which we may label as $[n_i]$. For any $t \in [n_i]$ the
outneighbourhood of $t$ in $H_i$ is $N^+_H(y'_{i,t}) \cap S'_{i+1,i}$ and the
inneighbourhood of $t$ in $H_i$ is $N^-(x'_{i,t}) \cap S'_{i-1,i}$. We also
include in $H_i$ all the edges from $S'_{i+1,i}$ to $S'_{i-1,i}$ that were present in $H$.
We claim that the underlying graph of each $H_i$ is $(d,18k\eta)$-super-regular.
Note that all three parts of $H_i$ have size $n_i > m/4k$, and therefore
any set containing at least $18k\eta n_i$ vertices from one of the parts
may be considered as a subset of $H$ of size at least $\eta m$.
Thus the regularity condition follows easily from the corresponding
$\eta$-regularity condition for $H$.
Next we need to check the degree condition in both
directions for each of the three directed bipartite graphs
$(Y'_i, S'_{i+1,i})$, $(S'_{i+1,i}, S'_{i-1,i})$, $(S'_{i-1,i}, X'_i)$.
Recall that the in and out neighbourhoods  of all vertices of $H$ in the set $S_{j,i}$ formerly had
size at least $(d-\eta)n_i - m^{2/3}$. Since $n_i > m/4k$ and we only swapped
$|B_{j,i}|<2\eta m+k<9k\eta n_i$ vertices, every vertex
still has at least $(d-10k\eta)n_i$ in and out neighbours in $S'_{j,i}$. Also, since vertices
in $S'_{i-1,i}$ are good they had at least $(d-2\eta)|X_i|$ outneighbours
in $X_i$. Again, we only removed at most
$|\cup_{j \ne i} B_{i,j}| \leq 2(2\eta m+k)<17k\eta n_i=17k\eta |X_i|$ vertices from $X_i$ to create $X'_i$.
Therefore, every vertex in $S'_{i-1,i}$ has at least
$(d-18k\eta)|X'_i|$ outneighbours in $X'_i$. Similar reasoning also shows that every vertex in
$S'_{i+1,i}$ has at least $(d-18k\eta)|Y'_i|$ inneighbours in $Y'_i$.
This establishes super-regularity. Now by Corollary \ref{tri-blowup} we can
perfectly cover each $H_i$ by vertex-disjoint cyclic triangles. This
translates into the required collection of  $k$-cycles in $H$. \qed

\medskip

Combining Theorem \ref{ck-blowup} with a technique similar to that used in
\cite{KO2, K} we can now prove Theorem \ref{ck}.

\medskip

\nib{Proof of Theorem \ref{ck}.}
Choose constants that satisfy
$1/n_0 \ll c \ll 1/M, 1/M' \ll \eps \ll d \ll \alpha \ll 1/k$.
Apply Lemma \ref{direg} to obtain a partition of the vertices of $G$ into
$V_0,V_1,\cdots,V_s$ and let $R'$ be the corresponding reduced digraph
with parameters $(\eps,d)$ on $[s]$. Let $R$ be the reduced oriented graph
obtained by applying Lemma \ref{red-orient} to $R'$. Then the minimum indegree
and outdegree in $R$ are at least $(1/2-c-3\eps-d)s > (1/2-2c)s$, so by
Theorem \ref{triangles} we can cover all but at most $\alpha s$ vertices
in $R$ by vertex-disjoint cyclic triangles $T_1,\cdots,T_t$.
(Although Theorem  \ref{triangles} allows us to cover all but at most $3$ vertices of $R$,
we only need this weaker bound which follows immediately from Theorem \ref{pip}).
Let $(V_{i,0},V_{i,1},V_{i,2})$ be the three clusters of the regular partition which correspond to the vertices of the
cyclic triangle $T_i$ in $R$. Since the oriented graph $R$ is a subgraph of the digraph $R'$ we have
that the edges of $G$ from $V_{i,j}$ to $V_{i,j+1}$ from an $\eps$-regular bipartite subgraph with density at least
$d$. By regularity, there are at most $2\eps |V_{i,j}|$ vertices $v \in V_{i,j}$ such that
$|N^+(v) \cap V_{i,j+1}|<(d-\eps)|V_{i,j+1}|$ or $|N^-(v) \cap V_{i,j-1}|<(d-\eps)|V_{i,j-1}|$.
We delete sets of size $\lfloor 2\eps |V_{i,j}| \rfloor$ from $V_{i,j}$ that contain all these vertices.
Thus we obtain triples $(U_{i,0},U_{i,1},U_{i,2})$ for $1 \le i \le t$ in
which each $U_{i,j}$ has the same size $m=\Omega(n/M)$,
and the edges between $U_{i,j}$ to $U_{i,j+1}$ are all
directed from $U_{i,j}$ to $U_{i,j+1}$.
Furthermore, since we deleted at most $2\eps$-proportion of each cluster, it is easy to see
that the underlying graph of edges between $U_{i,j}$ to $U_{i,j+1}$
forms a $(2\eps,d/2)$-super-regular pair.
Also, since $c \ll 1/M$, there is a constant $\delta \ll \eps$ such that
the minimum degree of all induced subgraphs $G[U_{i,j}]$ is at least
$|U_{i,j}|-cn \geq (1-\delta)|U_{i,j}|$.
Let $U_0$ denote the vertices that do not belong to any triple. Then $U_0$
contains the exceptional class $V_0$, which has size at most $\eps n$, the classes $V_i$ corresponding
to vertices of $R$ not covered by cyclic triangles, which have total size at most $\alpha n$, and the vertices
deleted to construct sets $U_{i,j}$, whose number is at most $2\eps \sum_i|V_i|=2\eps n$. Therefore
$|U_0| \leq \eps n+\alpha n+2\eps n<2\alpha n$.

Now we partition every set $U_{i,j}$ as $U'_{i,j} \cup U''_{i,j}$ by putting
each vertex randomly and  independently into either class with probability $1/2$.
Then, with high probability we have the following properties:
\begin{enumerate}
\item
$|U'_{i,j}|, |U''_{i,j}| = m/2 \pm m^{2/3}$ for every $i,j$,
\item
every vertex in $U_{i,j}$ has at least
$dm/5$ outneighbours in each of $U'_{i,j+1}$, $U''_{i,j+1}$
and inneighbours in each of $U'_{i,j-1}$, $U''_{i,j-1}$.
\item
for any vertex $x$, there are at least $(n/50)^{k-1}$ $k$-cycles
in which all vertices, except possibly $x$, are in $\cup_{i,j} U'_{i,j}$.
\end{enumerate}
The first two properties are simple applications of Chernoff bounds.
For the third property we use Azuma's inequality. Fix $x$ and
let $X$ be the random variable which counts the number of $k$-cycles
whose all vertices, except possibly $x$, are in $\cup_{i,j} U'_{i,j}$.
Since $|U_0| < 2\alpha n$ there are at most $2\alpha n^{k-1}$ such cycles containing $x$ and some vertex from $U_0$.
Therefore, by Lemma \ref{count-cycles} (proved in the next section) there are
at least $(n/10)^{k-1} - 2\alpha n^{k-1}$ $k$-cycles in which all vertices, except possibly $x$, are in $\cup_{i,j}
U'_{i,j}$. This implies that $\mb{E}X > 2^{-k+1}((1/10)^{k-1} - 2\alpha)n^{k-1}>(n/25)^{k-1}$.
Also, note that $X$ is a $c$-Lipchitz random variable with $c=n^{k-2}$. Indeed, there are most $n^{k-2}$ $k$-cycles
containing $x$ and any given vertex $v$. Hence
moving $v$ from $U'_{i,j}$ to $U''_{i,j}$ or vice versa can change the value of $X$ by at most $n^{k-2}$.
Now by Azuma's inequality (Theorem \ref{azuma}) we have
$$\mb{P}\big(X < (n/50)^{k-1}\big) < \mb{P}\big(|X-\mb{E}X| > (n/50)^{k-1}\big)
< 2e^{-(n/50)^{2(k-1)}/2n\cdot n^{2(k-2)}} < e^{-\Omega(n)}\ll 1/n\,. $$

Next we greedily cover the vertices in $U_0$ by disjoint $k$-cycles, so
that for each $x \in U_0$ we use $k$-cycle in which all other vertices
are in $\cup_{i,j} U'_{i,j}$. We do this in such a way to minimise the maximum
number of vertices used in any one of $U'_{i,j}$. When we come to cover
some $x \in U_0$, there are at least $(n/50)^{k-1}$ allowable $k$-cycles
(property 3 above). Of these,
at most $n^{k-2}|U_0| \leq 2 \alpha k n^{k-1}$ intersect a $k$-cycle that has already been
used to cover a vertex that came before $x$, and at most $\frac{2}{3}(n/50)^{k-1}$
intersect one of the heaviest (with respect to the number of vertices already used)
$50^{1-k} n/m$ classes $U'_{i,j}$
(since each $|U'_{i,j}|<2m/3$).
This means we can choose a $k$-cycle that is disjoint from those already chosen
and does not intersect one of the $50^{1-k} n/m$ heaviest classes,
and so the number of vertices used in any class will remain bounded by
$\frac{k|U_0|}{50^{1-k} n/m} < k50^k \alpha m$.

To finish the proof it is enough to show that when we restrict to the
uncovered vertices from each triple $(U_{i,0},U_{i,1},U_{i,2})$
we obtain a triple satisfying the hypotheses of Lemma \ref{ck-blowup}.
To see this recall that $|U_{i,j}|=m$ and $\alpha \ll 1$, so the number
of uncovered vertices in each class is at least
$m - k50^k \alpha m > 0.9m$. Super-regularity follows from
property 2 of the random partition, regularity of all pairs
$(U_{i,j}, U_{i,j+1})$ and the fact that we do not touch any vertex from
$(U''_{i,0},U''_{i,1},U''_{i,2})$. By Lemma \ref{ck-blowup}
we can cover all but at most $3k$ vertices in each triple by
disjoint $k$-cycles, so at most $C=3kt$ vertices remain uncovered. \qed

\section{Covering by prescribed cycles}

We have now assembled the two main ingredients for the proof of Theorem
\ref{almost-factor}, which we give in the first subsection of this section.
We have also done most of the preparation for the proof of Theorem
\ref{1-factor}: we will present a few more lemmas towards this end in
the second subsection, and then prove the theorem in the third subsection.

\subsection{Proof of Theorem \ref{almost-factor}.}

Choose $M$ so that Theorem \ref{long-cycles} applies with $\delta=1/10$
and then $c',C'$ so that Theorem \ref{ck} holds with parameters $c',C'$,
for all $k \le M$ and $n$ sufficiently large. Set $c=c'/2$, $C=MC'$.
Suppose that $n$ is sufficiently large, $G$ is an oriented graph on $n$ vertices
with minimum semidegree at least $(1/2-c)n$,
and $n_1,\cdots,n_t$ are numbers with $\sum_{i=1}^t n_i \le n-C$.
Let $N_k$ be the number of the $n_i$ equal to $k$, for $k \le M$, and
let $N_L = \sum_{n_i>M} n_i$. If there is any $k$ such that $N_k < cn/4M^2$
or if $N_L < cn/4$ we can greedily pack the appropriate cycles
using the previously mentioned Theorem 8 from \cite{KelKO2}, which says that an oriented graph
on $n$ vertices with minimum semidegree at least $3n/8$ contains cycles of all lengths
between $3$ and $n$. After that we will be left with minimum semidegree
at least $(1/2-c)n-\sum_{k=3}^M k (cn/4M^2)-cn/4 \geq (1/2-c)n-cn/4-cn/4 =(1/2-3c'/4)n$.
Thus we may reduce to the case when each $N_k$ for $k \le M$ is either $0$
or at least $cn/4M^2$ and $N_L$ is either $0$ or at least $cn/4$. Next we randomly
partition the remaining vertices, so that we allocate $kN_k+C'$ vertices
for the purpose of embedding $k$-cycles for each $k \le M$, and $N_L$ vertices
for the purpose of embedding all `long' cycles of length larger than $M$.
Lemma \ref{partition} implies that there is a choice of partition so
that each part has proportional semidegree at least $1/2-c'$, and then
Theorems \ref{ck} and \ref{long-cycles} allow us to embed the $k$-cycles
and the long cycles. This completes the proof. \qed

\subsection{Absorbing cycles}
When we have cycles of different lengths it is also useful to
consider the following kind of absorption. We say that a
cycle $F$ absorbs a path $P$ (disjoint from $F$) if
$F \cup P$ spans a (non-induced) cycle of length $|F|+|P|$. We present several lemmas
in this subsection that culminate in proving the
existence of a structure that is absorbing in this sense.

\begin{lemma} \label{join}
Suppose $0<c<10^{-4}$ and $G$ is an oriented graph on $n$ vertices
with minimum indegree and outdegree at least $(1/2-c)n$. Then $G$ has the following properties.

(1) Any $A \subset V(G)$ spans at least $e(A) \geq |A|(|A|/2-cn)$ edges.

(2) For any (not necessarily disjoint) subsets $S,T$ of $V(G)$
of size at least $(1/2-c)n$ there are at least $n^2/60$ directed edges from $S$ to $T$.

(3) For any (not necessarily disjoint) subsets $S,T$ of $V(G)$
of size at least $(1/2-c)n$ there are at least $10^{-5}n^3$ cyclic triangles
that contain an edge from $S$ to $T$.
\end{lemma}

\nib{Proof.}\, By deleting vertices if necessary we may assume that $|S|=|T|=(1/2-c)n$.

\noindent
(1) Since $|N(x)| \geq (1-2c)n$ for every vertex $x$ we obtain
\begin{eqnarray*}
e(A) &=& \sum_{x \in A} |N(x) \cap A|/2 \geq \sum_{x \in A} (|N(x)|+|A|-n)/2 \geq |A|((1-2c)n+|A|-n)/2\\
&=& |A|(|A|/2-cn) \,.
\end{eqnarray*}

\noindent
(2) Suppose first that $|S \cap T| > n/5$. Then, using the estimate from part (1), we get
$e(S,T) \ge e(S \cap T) \geq (n/5)(n/10-cn) > n^2/60$.
Otherwise $|\ov{S \cup T}| \leq n-(|S|+|T|-|S \cap T|) \leq n-\big(2(1/2-c)n-n/5\big) =(1/5+2c)n$. Therefore  we can
write
$$e(S,\ov{S}) = \sum_{x \in S} |N^+(x)| - e(S) > (1/2-c)n|S|-|S|^2/2=((1/2-c)n)^2/2$$
and
$$e(S,T) > e(S,\ov{S}) - |S||\ov{S \cup T}|> ((1/2-c)n)^2/2 - (1/2-c)n (1/5+2c)n > n^2/60.$$

\noindent
(3) From part (2) there are at least $n^2/60$ edges from $S$ to $T$.
By Lemma \ref{bad}, from any vertex $v\in S$ we have at most $(2a+4c)n$ outgoing edges which are $a$-bad. Taking
$a=1/300$ we obtain that at most $|S|n/100 \leq n^2/200$ of edges from $S$ to $T$ are
$1/300$-bad. Every $1/300$-good edge is contained in at least
$n/300$ cyclic triangles, each of which may be counted at most $3$ times,
so we get at least $(1/60 - 1/200)(1/900)n^3 > 10^{-5} n^3$ suitable
triangles. \qed

\begin{lemma} \label{count-cycles}
Suppose $0<c<10^{-4}$, $k \ge 3$ and $n$ is sufficiently large.
If $G$ is an oriented graph on $n$ vertices
with minimum indegree and outdegree at least $(1/2-c)n$ then
any vertex $x$ of $G$ belongs to at least $(n/10)^{k-1}$
$k$-cycles. More generally, if $t \ge 1$ and $k \ge t+2$ then any
path on $t$ vertices belongs to at least
$(n/10)^{k-t}$ $k$-cycles.
\end{lemma}

\nib{Proof.} To construct a $k$-cycle through $x$ we start
by greedily picking a path of $k-2$ vertices starting at $x$.
When $k=3$ this is just the point $x$. For $k \ge 4$, note that
by the outdegree condition we have at least $(1/2-c)n-k$ choices at every step,
so this gives at least $\prod_{i=0}^{k-4} ((1/2-c)n-k)$ such paths.
Given a path $P$ of length $k-2$, from $x$ to some final point $y$, we may complete
$P$ to a $k$-cycle by choosing an edge from $N^+(y)$
to $N^-(x)$ which does not use any vertex of $P$. Clearly there are at most $kn$ edges incident to the
vertices on the path $P$. Hence, by Lemma \ref{join},
there are at least $n^2/60 - kn$ edges from $N^+(y)$
to $N^-(x)$ disjoint from $P$. Altogether
we get at least $(n^2/60 - kn)\prod_{i=0}^{k-4} ((1/2-c)n-k)> (n/10)^{k-1}$ cycles. The estimate
for the number of $k$-cycles containing a given path of length $t$ can be obtained similarly. \qed

\begin{lemma} \label{ck-cl}
Suppose $0<c<10^{-4}$, $k \ge 3$, $\ell \ge k+3$ and $n$ is sufficiently large.
If $G$ is an oriented graph on $n$ vertices
with minimum indegree and outdegree at least $(1/2-c)n$
and $P$ is any path on $k$ vertices in $G$ then there
are at least $(n/100)^{\ell-k}$ cycles $C$ of length $\ell-k$ so that
$P \cup C$ spans a (non-induced) cycle of length $\ell$.
\end{lemma}

\nib{Proof.}
Let $S$ be the outneighbourhood of the last vertex of $P$
and $T$ the inneighbourhood of the first vertex of $P$.
Suppose first that $\ell > k+3$. By part (2) of Lemma \ref{join}
there are at least $n^2/60$ edges $xy$ with $x \in T$ and $y \in S$.
Also, by Lemma \ref{count-cycles} each
such $xy$ is contained in at least $(n/10)^{\ell-k-2}$ cycles of length $\ell-k$.
Altogether this gives at least $(n^2/60)(n/10)^{\ell-k-2}$ cycles of length $\ell-k$.
Since  at most $kn^{\ell-k-1}$ of these cycles intersect the
path  $P$, there are at least
$ (n^2/60)(n/10)^{\ell-k-2} -kn^{\ell-k-1}> (n/10)^{\ell-k}$
cycles of length $\ell-k$ containing an edge from $T$ to $S$ and
disjoint from $P$. Clearly, each such cycle together with $P$ spans a cycle of length $\ell$.

Now suppose that $\ell=k+3$. By assertion (3) of Lemma \ref{join}  there are
at least $10^{-5}n^3$ cyclic triangles that contain an edge from $T$ to $S$. At most
$kn^2$ of these triangles use a point from $P$, so at least
$10^{-5}n^3-kn^2>10^{-6}n^3 = (n/100)^{\ell-k}$ are disjoint from $P$.
These triangles together with $P$ span cycles of length $\ell$.
\qed

\medskip

Now by Lemmas \ref{count-cycles} and \ref{ck-cl},
the same argument that we used in Lemma \ref{absorb-tri},
using simply Chernoff bounds rather than the Kim-Vu inequality,
leads to the following lemma.

\begin{lemma} \label{absorb-ck-cl}
For any $k \ge 3$ and $M \ge k+1$
there is some $c>0$ and number $n_0$ such that
if $G$ is an oriented graph on $n>n_0$ vertices
with minimum indegree and outdegree at least $(1/2-c)n$
then the following holds.
Suppose we form a collection of vertex-disjoint $k$-cycles
$\cal C$ by choosing each $k$-cycle independently with some probability
$p$
and deleting any pair of $k$-cycles that intersect. Write $a_k n^k$ for the number
of $k$-cycles in $G$ (where $a_k>k^{-1}10^{-k+1}$ by Lemma
\ref{count-cycles}).
If $\frac{\log n}{n^k} \ll p \ll 1/n^{k-1}$ then with high probability
we have $|{\cal C}| = m \sim a_k p n^k$ and for any path $P$ on $\ell-k$
vertices
with $k+1 \le \ell \le M$ there are at least $200^{-k} m$
absorbing $k$-cycles for $P$ in $\cal C$.
\end{lemma}

\subsection{Proof of Theorem \ref{1-factor}.}

Choose constants with the hierarchy
$M \ll c^{-1} \ll C \ll T \ll n_0$.
By assumption at least $T$ of the $n_i$
lie between $k+1$ and $M$. We may relabel so that $k+1 \le n_i \le M$
for $1 \le i \le T$. Next by Lemma \ref{absorb-ck-cl} we choose a collection $\cal C$ of
$|{\cal C}|=T \log n + T$ vertex-disjoint $k$-cycles such that
for any path $P$ in $G$ on $\ell-k$ vertices
with $k+1 \le \ell \le M$ there are at least $200^{-k} m$
absorbing $k$-cycles for $P$ in $\cal C$. Let $G'$ be the restriction of $G$ to the vertices not covered by
cycles in $\cal C$. This is a tournament
on $n' = n - kT(\log n +1)$ vertices with minimum semidegree
at least $(1/2-c)n- kT(\log n +1)>(1/2-2c)n'$.
Let $n'_1,\cdots,n'_{t'}$ be the sequence obtained from $n_1,\cdots,n_t$
by removing $n_1,\cdots,n_T$ and $T \log n$ occurrences of $k$.
Note that $\sum_j n'_j=n'-\sum_{i=1}^T (n_i-k) < n'-C$. Therefore, we can
apply Theorem \ref{almost-factor} and find a packing of cycles in $G'$ of length
$n'_1,\cdots,n'_{t'}$ covering all but a set $U$ of $\sum_{i=1}^T (n_i-k)$
vertices.
Note that $G'$ restricted to $U$ is a tournament, and so contains a
Hamilton path. We can partition this path into $T$ paths with $n_1-k,\cdots,n_T-k$
vertices. Finally, we can apply the absorbing property of $\cal C$ to repeatedly
combine a path on $n_i-k$ leftover vertices with a $k$-cycle in $\cal C$
to form an $n_i$-cycle, for $1 \le i \le T$. This completes the proof. \qed

\section{Concluding remarks}

In \cite{C}, Cuckler raises the question of counting perfect
packings of $k$-cycles in regular tournaments (when $k$ does not divide
$n$ a `perfect' packing is defined to have size $\lfloor n/k \rfloor$).
He conjectures  that for odd $k$ the number of such packings is
$n!^{(k-1)/k}(2+o(1))^{-n}$, which is asymptotically the number of
perfect $k$-cycle packings which one expects to have in a random tournament.
Somewhat surprisingly,  he shows that this is no longer true if
$k$ is even. In the same paper Cuckler also gives this estimate for counting
triangle packings of size $n/3 - o(n/\log n)$ in regular tournaments.
Our proof of Theorem
\ref{triangles} can be used to show that any oriented graph $G$
of order $n$ with semidegree $(1/2+o(1))n$ has
$n!^{2/3}(2+o(1))^{-n}$ triangle packings covering all but at most $3$
vertices (which are `perfect' when $n$ is not divisible by $3$).
The upper bound on the number of packings follows simply from the fact that the number of cyclic triangles
in any tournament of order $n$, and hence also in any oriented graph, is at most
$(1+o(1))n^3/24$. Constructing a packing by choosing one triangle every time we see that
for the $i$-th triangle we have at most
$(1+o(1))(n-3(i-1))^3/24$ choices. Dividing by the number of different orderings of the same
packing, which is at least $(n/3-1)!$, we see that there are at most
$$ \frac{(1+o(1))^n}{(n/3-1)!}\prod_{i=1}^{n/3}\frac{(n-3(i-1))^3}{24} = \frac{(1+o(1))^n}{(n/3-1)!}
2^{-n}3^{-n/3}n!=n!^{2/3}(2+o(1))^{-n}$$
different packing of cyclic triangles covering all but at most $3$ vertices.

Next we present a sketch proof for the lower bound.
Suppose that the nibble consists of $m$ `bites' of size
$b_1,\cdots,b_m$ with each $b_i=o(n)$ and
$\sum b_i = n/3-o(n)$. From the analysis of the nibble in \cite{AS} it is easy to see that
at iteration $i$ we have an
oriented graph on $n-\sum_{j<i}3b_j$ vertices which has
$(1+o(1)) \frac{(n-\sum_{j<i}3b_j)^3}{24}$ cyclic triangles; then we pick $b_i$ such
triangles uniformly at random. Therefore, the number of ordered
choices of $n-o(n)$ vertex-disjoint cyclic triangles in $G$ is at least
$(1+o(1))^n \prod_{i=1}^m \left((n-\sum_{j<i}3b_j)^3/24\right)^{b_i}$.
Since for $b=o(a)$ we have
$$(a^3/24)^b =(1+o(1))^b24^{-b}a(a-1)\cdots(a-3b+1),$$
we conclude that
the number of ordered choices of $n-o(n)$ vertex-disjoint cyclic triangles in $G$ is
at least
$$(1+o(1))^n24^{-\sum_ib_i}n!=
(1+o(1))^n24^{-n/3}n!\,.$$
We can use the absorption argument to convert each of these families of
$n-o(n)$ vertex-disjoint cyclic triangles into triangle packings
that cover all but at most $3$ vertices. Note that since the absorption process involves only
$o(n)$ vertices, each such packing is obtained
in at most $\binom{n/3}{o(n)}=e^{o(n)}$ ways. Note also that for each packing its triangles can
have at most
$(n/3)!$ different orders. Dividing by
this number we obtain that the number of
triangle packings covering all but at most $3$ vertices in $G$ is at least
$(1+o(1))^n\, e^{-o(n)}\, 24^{-n/3} n! \big/ (n/3)! = (1/2+o(1))^n n!^{2/3}$, as required.

As discussed in the introduction, it is also natural to consider questions which
we study in this paper for digraphs, rather than oriented graphs.
It appears that such digraph problems tend to be quite closely connected
to known extremal results about graphs, and can therefore be answered using these results.
Here we mention two illustrative examples:
\begin{enumerate}
\item
A digraph $G$ on $n$ vertices with minimum total degree (this is a minimum
of $d^+(x)+d^-(x)$ over all $x \in G$) at least
$(4/3+c)n$ has a perfect packing of transitive triangles, when $n>n_0(c)$ is
divisible by $3$ and sufficiently large.
To see this consider a random ordering $<$ of the vertices of $G$ and let $G'$ be the graph whose
edges are all those edges $ij$ of $G$ with $i<j$. It is easy to check that with high probability
every vertex of $G'$ has degree at least $(2/3+c/4)n$.
Now the Corr\'adi-Hajnal theorem (mentioned in the introduction) gives a perfect triangle packing in $G'$, which by definition
corresponds to a packing of digraph $G$ by transitive triangles.

\item
A digraph $G$ on $n$ vertices with minimum semidegree at least $(2/3+c)n$
has a perfect packing of cyclic triangles, when $n>n_0(c)$ is
divisible by $3$ and sufficiently large.  To see this, randomly partition
the vertices of $G$ into three disjoint sets $V_0 \cup V_1 \cup V_2$ with $|V_i|=n/3$, $i=0,1,2$.
Let $G'$ be the $3$-partite graph whose edges are the edges of $G$
which go from $V_i$ to $V_{i+1}$ (addition mod $3$).
Again, it is easy to check using large deviation inequalities that every vertex of $G'$ in $V_i$
has at least $(2/3+c/2)|V_j|$ neighbours in each $V_j$, $j \ne i$.
A result of Johansson \cite{J} implies that $G'$ has a perfect triangle
packing, which by definition gives a cyclic triangle packing in $G$.
\end{enumerate}
Moreover, one can easily show that both these results are asymptotically best
possible by taking corresponding construction for graphs and replacing each edge
by two directed edges with opposite directions.

Note that in Theorem \ref{triangles} we find a triangle packing which covers all but
three vertices. What happens if we slightly relax this requirement?
In particular, what minimum semidegree condition in an oriented graph of order $n$
will give a cyclic triangle packing that covers all but at most $o(n)$ vertices?
We have no good conjecture for this problem.
The following construction shows that the semidegree should be at least $(4/9-o(n))$.
Suppose $\eps>0$ and divide a set of $n$ vertices into
three sets $V_0,V_1,V_2$ with sizes $|V_0|=(1/3-\eps/2)n$, $|V_1|=n/3$,
$|V_2|=(1/3+\eps/2)n$. Define an oriented graph $G$ as follows.
Between the classes we take all possible edges and orient them from $V_i$ to $V_{i+1}$ (addition mod $3$).
Inside each class we place an oriented graph that has
no cyclic triangle with minimum indegree and outdegree as large
as possible. For example, a circulant construction gives
such an oriented graph on $m$ vertices with every indegree
and outdegree equal to $\lceil m/3 \rceil - 1$, and this
cannot be improved if the Caccetta-Haggkvist conjecture \cite{CaH} is true.
Then $G$ has minimum indegree and outdegree at
least $(4/9 - \eps)n$. Since any cyclic triangle must use
one vertex from each class, any collection of
vertex-disjoint cyclic triangles leaves at least $\eps n$
vertices uncovered.

Another interesting variation is to consider what minimum semidegree
condition is needed to find certain structures in a tournament. Here one would
expect a smaller value than that needed for oriented graphs. For example, using the fact that every strongly connected
tournament is Hamiltonian (Camion's theorem), it is not hard to see that a tournament with minimum
semidegree at least $n/4$ contains a Hamilton cycle. On the other hand, we recall that for oriented graphs
ones need semidegree at least $3n/8$ to get the same conclusion.
It would be interesting to find what minimum semidegree in a tournament of order $n$
will give a cyclic triangle packing that covers all but at most $o(n)$ vertices.
Here it is easy to show a lower bound of $n/3$ (note it is again smaller than that
obtained above for oriented graphs).

\end{document}